\documentclass[a4paper,10pt]{article}

\usepackage{amssymb,amsmath,graphics}

\newenvironment{dwd}{\par\noindent{\bf Proof.}}{\par\rightline{$\blacksquare$}}

\newtheorem{theo}{Theorem}
\newtheorem{prop}{Proposition}  
\newtheorem{coro}{Corollary}

\newtheorem{lema}{Lemma}
\newtheorem{defi}{Definition}
\def\be#1\ee{\begin{equation}#1\end{equation}}
\newcommand{\ba}{\begin{eqnarray} }
\newcommand{\ea}{\end{eqnarray} }
\def\bt#1\et{\begin{theo}#1\end{theo}}
\def\bl#1\el{\begin{lema}#1\end{lema}}
\def\bp#1\ep{\begin{prop}#1\end{prop}}
\def\bd#1\ed{\begin{defi}#1\end{defi}}

\def\ccA{{\cal A}}
\def\ccB{{\cal B}}

\def\va{\varepsilon}
\def\ra{\rightarrow}

\def\E{\mathbf{E}}
\def\P{\mathbf{P}}

\def\N{{\mathbb N}}

\def\Z{{\mathbb Z}}

\def\ls{\leqslant}
\def\gs{\geqslant}

\begin{document}

\title{\bf Some remarks on the Oleszkiewicz problem
}
\author{Witold Bednorz
\footnote{Partially supported Research partially 
supported by the NCN grant DEC-
2012/05/B/ST1/00412..}
\footnote{Institute of Mathematics, University of Warsaw, Banacha 2, 02-097 Warszawa, Poland}}

\maketitle
\begin{abstract}
In this paper we study the question how to easily verify that the expectation
of the supremum of a one canonical Bernoulli process dominates the same quantity for another process of this type. In the setting of Gaussian canonical processes it is known 
that such a comparison holds for contractions in the Euclidean distance. We do state and prove a similar result for Bernoulli processes. 
In particular we get a partial answer to the Oleszkiewicz conjecture about the comparability of weak and strong moments for type Bernoulli series in a Banach space.  
\end{abstract}

\section{Introduction}

Suppose that $T\subset \ell^2$ and let $X_t$, $t\in T$ be defined as
$$
X_t=\sum_{i\gs 1}t_i\va_i,\;\;t=(t_i)_{i\gs 1},
$$
where $\va_i$ are independent random signs, i.e. $\P(\va_i=\pm 1)=\frac{1}{2}$.
Let 
$$
b(T)=\E\sup_{t\in T}X_t=\E\sup_{t\in T}\sum_{i\gs 1}t_i\va_i.
$$
Suppose that $f:T\ra \ell^2$, we aim to formulate a simple conditions on such a function $f$
which guarantees that $b(f(T))\ls Kb(T)$. Note that the similar question can be posed in the setting of Gaussian canonical processes. Namely assume that $T\subset \ell^2$ and let $G_t$, $t\in T$
be defined as
$$
G_t=\sum_{i\gs 1}t_i g_i,\;\;t=(t_i)_{i\gs 1},
$$ 
where $g_i$ are independent standard normal variables. Let 
$$
g(T)=\E\sup_{t\in T}G_t=\E\sup_{t\in T}\sum_{i\gs 1}t_ig_i.
$$ 
Then it is well-known  that the following holds.
\bt\label{thm:4}
Suppose that $\pi:T\ra \ell^2$ satisfies
$$
\|f(t)-f(s)\|_2\ls C\|t-s\|_2,\;\;\mbox{for all}\;\;s,t\in T
$$
then $g(f(T))\ls Kg(T)$, where $K$ is a universal constant.
\et
In the case of Bernoulli processes the answer must be more complicated since
the distribution of $X_t$ much more depends on the structure of $t=(t_i)_{i\gs 1}$, namely it is
known by \cite{Hit} that for any $p\in \N$, $p\gs 1$
\be\label{ole0}
\|X_t\|_p\sim \sum^p_{i=1}|t^{\ast}_i|+\sqrt{p}(\sum_{i>p}|t^{\ast}_i|^2)^{\frac{1}{2}},
\ee
where $(t^{\ast}_i)_{i\gs 1}$ is the rearrangement of $(t_i)_{i\in I}$ such that
$|t^{\ast}_1|\gs |t^{\ast}_2|\gs \ldots$. It is obvious that
$$
\P(|X_t|\gs e\|X_t\|_p)\ls e^{-p},\;\;p\gs 1
$$
on the other hand there exists constant $\kappa$
$$
\P(|X_t|\gs \kappa^{-1}\|X_t\|_p)\gs \min\{\kappa^{-1},e^{-p}\},\;\;p\gs 1.
$$
Therefore the simplest condition one could formulate when comparing 
$(X_{f(t)})_{t\in T}$ and $(X_t)_{t\in T}$ is
\be\label{ole1}
\|X_{f(t)}-X_{f(s)}\|_p\ls C\|X_{t}-X_{s}\|_p,\;\;s,t\in T,\;\;p\gs 1.
\ee
Up to now there is no counterexample to the conjecture that (\ref{ole1}) implies $b(f(T))\ls Kb(T)$,
where $K$ is a universal constant but it seems that there is no appropriate approach
to establish such a result. A weak support to the conjecture
is the following chaining argument. We say that $\ccA=(\ccA_n)_{n\gs 0}$ is an admissible
sequence of partitions of $T$ if $\ccA_{0}=\{T\}$ and $|\ccA_n|\ls N_n=2^{2^n}$ for $n\gs 1$
(usually there is required also that these partitions are nested).
Let $\pi_n(A)$ be a given point in $A\in \ccA_n$ and $\pi_n(t)=\pi_n(A_n(t))$, where
$t\in A_n(t)\in \ccA_n$. Let
$$
\gamma_X(T)=\inf_{\ccA}\sup_{t\in T}\sum_{n\gs 1}\|X_{\pi_n(t)}-X_{\pi_{n-1}(t)}\|_{2^n}\|,
$$
where the infimum is taken over all admissible sequences of partitions.
It is well-known that $b(T)\ls K\gamma_X(T)$, where $K$ is a universal constant.
Consequently (\ref{ole1}) implies the following result.
\bt\label{thm:0}
Suppose that (\ref{ole1}) holds then
$$
b(f(T))\ls K\gamma_X(f(T))\ls KC\gamma_X(T),
$$
where $K$ is a universal constant.
\et
However in general $\gamma_X(T)$ is not comparable with $b(T)$.
The question how to characterize $b(T)$ up to a universal constant was for a long time 
an open question finally solved in \cite{Bed1} but still what lacks is some better understanding
of the result. The main idea how to bound $b(T)$ from below is that there should exists 
a decomposition of $T$ into $T_1,T_2\subset \ell^2$ so that $T_1+T_2=\{s+t\in \ell^2:\; s,t\in T\}$ covers set $T$. Usually such a decomposition is formulated in the language of existence of $\pi:T\ra \ell^2$ which is used to define $T_1=\{t-\pi(t):t\in T\}$ and $T_2=\{\pi(t):t\in T\}$.
We have the following consequence of the Bernoulli Theorem \cite{Bed1}.
\bt\label{thm:3}
There exists a function $\pi:T\ra \ell^2$ such that
$$
K^{-1}(\gamma_X(T_1)+\gamma_X(T_2))\ls b(T)\ls K(\gamma_X(T_1)+\gamma_X(T_2)),
$$
where $K$ is a universal constant, $T_1=\{t-\pi(t):t\in T\}$ and $T_2=\{\pi(t):t\in T\}$.
\et
\begin{dwd}
By the main result of \cite{Bed1} we get the existence of $\pi:T\ra \ell^2$ and consequently $T_1$ and $T_2$ such that
\be\label{ole3}
b(T)\gs L^{-1}(\sup_{t\in T_1}\|t\|_1+\gamma_G(T_2)),
\ee
where $\|t\|_1=\sum_{i\gs 1}|t_i|$ and 
$$
\gamma_G(T)=\inf_{\ccA}\sup_{t\in T}\sum_{n\gs 1}\|G_{\pi_n(t)}-G_{\pi_{n-1}(t)}\|_{2^n},
$$
where the infimum is taken over all admissible sequences of partitions.
It is obvious that $\gamma_G(T_2)\gs \gamma_{X}(T_2)$. On the other hand
it has to be noticed that $\sup_{F\subset T_1}\gamma_X(F)=\gamma_X(T_1)$
where the supremum is over all finite subsets in $T_1$. 
Then it suffices to create trivial partition sequence $\{F\}=\ccA_0=\ccA_1=\ldots=\ccA_{M-1}$
and $\ccA_M=F$ where $N_M\gs |F|>N_{M-1} $. Let $t_0=\pi_0(F)$ and observe that for any $t\in F$
$$
\|X_t-X_{t_0}\|_{2^M}\ls \|t-t_0\|_1\ls 2\sup_{t\in F}\|t\|_1 
$$
and hence
$$
\gamma_X(T_1)\ls 2\sup_{t\in T_1}\|t\|_1.
$$ 
Consequently
$$
b(T)\gs K^{-1}(\gamma_X(T_1)+\gamma_X(T_2)).
$$
On the other hand it is trivial to get
$$
b(T)\ls b(T_1)+b(T_2)\ls K(\gamma_X(T_1)+\gamma_X(T_2)).
$$
\end{dwd}
In general we do claim that for any canonical process
of the from
$$
Y_t=\sum_{i\gs 1}t_i\xi_i,\;\;t=(t_i)_{i\gs 1},
$$
where $t\in T\subset\ell^2$ and $\xi_i$ are independent symmetric identically distributed random variables
such that $\E \xi_i^2=1$ and whose distribution has log-concave tails 
there exists a function $\pi:T\ra\ell^2$ such that
$$
K^{-1}(\gamma_Y(T_1)+\gamma_Y(T_2))\ls\E\sup_{t\in T}Y_t\ls K(\gamma_Y(T_1)+\gamma_Y(T_2)),
$$
where $T_1=\{t-\pi(t)\in \ell^2:\;t\in T\}$, $T_2=\{\pi(t)\in \ell^2:\;t\in T\}$ and
$$
\gamma_Y(T)=\inf_{\ccA}\sup_{t\in T}\sum_{n\gs 1}\|Y_{\pi_n(t)}-Y_{\pi_{n-1}(t)}\|_{2^n},
$$
where the infimum is taken over all admissible sequences of partitions. That is why
the question we treat in this paper can be studied in much extent.
Unfortunately there is no canonical way to find the function $\pi:T\ra\ell^2$ and  this 
makes the problem difficult when we try to transport the construction by the map $f$. 
Therefore we prove a weaker form of our conjecture. Note that (\ref{ole1}) means that
$$
\|\sum_{i\gs 1}|f(t)_i-f(s)_i|\va_i\|_p\ls C\|\sum_{i\gs 1}|t_i-s_i|\va_i\|_p,\;\;s,t\in T,\;\;p\gs 1.
$$
Our condition will be based on the truncation of coefficients in the above inequality, namely
we prove the following result.
\bt\label{thm:1}
Suppose that for all $s,t\in T$, $p\gs 1$ and $r>0$
\be\label{ole2}
\|\sum_{i\gs 1}(|f(t)_i-f(s)_i|\wedge r)\va_i\|_p\ls C(rp+\|\sum_{i\gs 1}(|t_i-s_i|\wedge r)\va_i\|_p).
\ee
\et
The result is stronger then what can be derived from the Bernoulli comparison - Theorem 2.1 in \cite{Tal2} (also Theorem 5.3.6 in \cite{Tal1})  where it is assumed that 
for all $t\in T$ and $i\gs 1$, $|f(t)_i-f(s)_i|\ls |t_i-s_i|$. In fact we have a simple corollary of Theorem \ref{thm:1}.
\begin{coro}
Suppose that for any $p\in \N$ and $C\gs 1$
\be\label{mela}
\inf_{I:|I|\ls Cp}(\sum_{i\not\in I}|f(t)_i-f(s)_i|^2)^{1/2}\ls C\inf_{I:|I|=p}(\sum_{i\not\in I}|t_i-s_i|^2)^{1/2}.
\ee
Then $b(f(T))\ls Kb(T)$, where $K$ is a universal constant.
\end{coro}
\begin{dwd}
Indeed it suffices to use (\ref{ole0}) to get that (\ref{mela}) implies (\ref{ole2}) and then apply Theorem \ref{thm:1}.
\end{dwd}
One of the questions which is related to the above analysis is the comparison of weak and strong moments for
type Bernoulli series. We will describe the idea in the last section.

\section{Proof of the main result}

We prove in this section Theorem \ref{thm:1}.   
\begin{dwd}
The main step in the proof of Bernoulli theorem is to show the existence of a suitable admissible
sequence partition. Consequently if $b(T)<\infty$ and say $0\in T$ then it is possible to define
nested partitions $\ccA_n$ of $T$ such that $|\ccA_n|\ls N_n$. Moreover for each $A\in \ccA_n$ it possible to find $j_n(A)\in \Z$ and $\pi_n(A)\in T$ (we use the 
notation $j_n(t)=j_n(A_n(t))$ and $\pi_n(t)=\pi_n(A_n(t))$) which satisfies
the following assumptions
\begin{enumerate}
\item $\|t-s\|_2\ls \sqrt{M}r^{-j_0(T)}$, for $s,t\in T$;
\item if $n\gs 1$, $\ccA_n\ni A\subset A'\in \ccA_{n-1}$ then 
\begin{enumerate}
\item either $j_n(A)=j_{n-1}(A')$ and $\pi_n(A)=\pi_{n-1}(A')$ 
\item or $j_n(A)>j_{n-1}(A')$, $\pi_n(A)\in A'$ and
$$
\sum_{i\in I_n(A)}\min\{|t_i-\pi_n(A)_i|^2,r^{-2j_n(A)}\}\ls M2^nr^{-2j_n(A)}
$$ 
where for any $t\in A$
$$
I_n(A)=I_n(t)=\{i\gs 1:\;|\pi_{k+1}(t)_i-\pi_k(t)_i|\ls r^{-j_k(t)}\;\;\mbox{for}\;\;0\ls k\ls n-1\}
$$
\end{enumerate}
\item Moreover numbers $j_n(A)$, $A\in \ccA_n$, $n\gs 0$ satisfy
\be\label{ole4}
\sup_{t\in T}\sum^{\infty}_{n=0}2^nr^{-j_n(t)}\ls Lb(T).
\ee
\end{enumerate}
As it is proved in Theorem 3.1 in \cite{Bed1} the existence of such a construction implies the existence of a decomposition $T_1,T_2\subset \ell^2$, $T_1+T_2\supset T$ such that
$$
\sup_{t^1\in T_1}\|t^1\|_1\ls LM\sup_{t\in T}\sum^{\infty}_{n=0}2^nr^{-j_n(t)}\;\;\mbox{and}\;\;
\gamma_G(T_2)\ls L\sqrt{M}\sup_{t\in T}\sum^{\infty}_{n=0}2^nr^{-j_n(t)}.
$$    
In this why we get (\ref{ole3}). Now if we have mapping $f:T\ra f(T)\subset \ell^2$
we can preserve a lot of properties of the construction of $\ccA_n$, $\pi_n(A)$ and $j_n(A)$.
Namely let $\ccB_{n}$ consists of $f(A)$, $A\in \ccA_n$. Obviously partitions $\ccB_n$ are admissible, nested
and $\ccB_0=\{f(T)\}$. Then we define for each $n\gs 0$ and $A\in \ccB_n$
$$
\pi_n(f(A))=f(\pi_n(A))\;\;\mbox{and}\;\;j_n(f(A))=j_n(A).
$$ 
We have to verify all the assumptions from Theorem 3.1 in \cite{Bed1}.
To this aim we need our main condition (\ref{ole2}). Let $C_0\gs C$ be suitably large constant.
First it is obvious for large enough $r=\sqrt{M}r^{-j_0(T)}$ and $p=2$ that (\ref{ole2}) implies
$$
\|f(t)-f(s)\|_2\ls 4C(\sqrt{M}r^{-j_0(T)}+\|t-s\|_2)\ls 8C_0\sqrt{M}r^{-j_0(T)}.
$$
Then if $f(A)\in\ccB_n$ and $f(A)\subset f(A')\in \ccB_{n-1}$ then either
$$
j_n(f(A))=j_n(A)=j_{n-1}(A')=j_{n-1}(f(A'))
$$
and
$$
\pi_n(f(A))=f(\pi_n(A))=
f(\pi_{n-1}(A'))=\pi_{n-1}(f(A'))
$$
or $j_n(f(A))=j_n(A)>j_{n-1}(A')=j_{n-1}(f(A'))$. In this case we have
$\pi_n(f(A))=f(\pi_n(A))\in f(A')$ and it suffices to show that
\be\label{ole5}
\sum_{i\in I_n(f(A))}\min\{|f(t)_i-f(\pi_n(A))_i|^2,r^{-2j_n(f(A))}\}\ls 64C^2_0M2^nr^{-2j_n(f(A))}
\ee
where 
$$
I_n(f(A))=I_n(f(t))=\{i\gs 1: |f(\pi_{k+1}(t))_i-f(\pi_k(t))_i|\ls r^{-j_k(f(t))}\;\mbox{for}\;0\ls k\ls n-1\}.
$$
To establish (\ref{ole5}) we first observe that 
\begin{align*}
& \sum_{i\in I_n(f(A))}\min\{|f(t)_i-f(\pi_n(A))_i|^2,r^{-2j_n(f(A))}\}\\
&\ls \sum_{i\gs 1}\min\{|f(t)_i-f(\pi_n(A))_i|^2,r^{-2j_n(f(A))}\}
\end{align*}
Now due to (\ref{ole1}) and (\ref{ole2}) we deuce that
\begin{align*}
&  \sum_{i\gs 1}\min\{|f(t)_i-\pi_n(A)_i|^2,r^{-2j_n(f(A))}\} \\
& \ls C_1(2^nr^{-2j_n(f(A))}+\|\sum_{i\gs 1}\min\{|f(t)_i-f(\pi_n(A))_i|^2,r^{-j_n(f(A))}\}\va_i\|^2_{2^n})\\
& \ls C_2(2^nr^{-2j_n(f(A))}+\|\sum_{i\gs 1}\min\{|t_i-\pi_n(A)_i|^2,r^{-j_n(f(A))}\}\va_i\|^2_{2^n})\\
& \ls C_3(2^nr^{-2j_n(A)}+\sum_{i\gs 1}\min\{|t_i-\pi_n(A)_i|^2,r^{-2j_n(A)}\}), 
 \end{align*}
where in the last line we have used that $j_n(f(A))=j_n(A)$. It remains to
observe that $|I_n(A)^c|\ls C_42^n$. Let $t\in A$ if $\pi_{k+1}(t)\neq \pi_k(t)$ then
$j_{k+1}(t)>j_k(t)$ and hence $\pi_{k+1}(t)\in A_k(t)$. Therefore there exists $l\in \{0,1,\ldots,k\}$
such that
$$
j_k(t)=j_{k-l}(t)>j_{n-l-1}(t)
$$  
and hence $\pi_{k+1}(t)\in A_{k-l}(t)$ and $\pi_k(t)=\pi_{k-l}(t)$, $j_k(t)=j_{k-l}(t)$ so by the
construction of $(\ccA_n)_{n\gs 0}$
\begin{align*}
& \sum_{i\in I_{k-l}(t)}\min\{(\pi_{k+1}(t)-\pi_{k-l}(t))^2,r^{-2j_{k-l}(t)}\}\\
&=\sum_{i\in I_{k-l}(t)}\min\{(\pi_{k+1}(t)-\pi_{k}(t))^2,r^{-2j_k(t)} \}\ls M2^kr^{-2j_k(t)}.
\end{align*}
Consequently
$$
|\{i\in I_{k-l}(t):\;|\pi_{k+1}(t)_i-\pi_k(t)_i|>r^{-j_k(t)}\}|\ls M2^k. 
$$
Therefore  $|I^c_n(A)|\ls M\sum^n_{k=0}2^{k}\ls M2^{n+1}$ and hence
\begin{align*}
& \sum_{i\gs 1}\min\{|f(t)_i-\pi_n(A)_i|^2,r^{-2j_n(A)}\}\ls \\
&\ls C_4(2^nr^{-2j_n(A)}+\sum_{i\in I_n(A)}\min\{|t_i-\pi_n(A)_i|^2,r^{-2j_n(A)}\})\ls 
64C_0^2 M^22^n r^{-2j_n(A)},
\end{align*}
for suitably large constant $C_0$. All the assumptions required in Theorem 3.1 in \cite{Bed1} are satisfied for $(\ccB_n)_{n\gs 0}$.
Therefore there exists a decomposition $S_1,S_2\subset \ell^2$ such that $S_1+S_2\supset f(T)$
and
$$
\sup_{s\in S_1}\|s\|_1\ls 8C_0LM\sup_{t\in f(T)}\sum_{n\gs 0}2^nr^{-j_n(t)},\;\;
\gamma_G(S_2)\ls 8C_0L\sqrt{M}\sup_{t\in f(T)}\sum_{n\gs 0}2^nr^{-j_n(t)}. 
$$
Since $j_n(f(t))=j_n(t)$ we get by (\ref{ole4}) that
$$
\sup_{t\in f(T)}\sum_{n\gs 0}2^nr^{-j_n(t)}\ls Lb(T).
$$
It implies that
$$
b(f(T))\ls b(S_1)+b(S_2)\ls Kb(T),
$$
for a universal constant $K$ and ends the proof.
\end{dwd}

\section{The Oleszkiewicz problem}

 Our study can be applied to the question posed by Krzysztof Oleszkiewicz that concerns comparability of weak and strong moments for type Bernoulli series in a Banach space.  
 Let $x_i$, $y_i$, $i\gs 1$ be vectors in  a Banach space $(B,\|\; \|)$. 
 Suppose that for all $x^{\ast}\in B^{\ast}$ and $u\gs 0$
 \be\label{Ole1}
 \P(|\sum_{i\gs 1}x^{\ast}(x_i)\va_i|>u)\ls C\P(|\sum_{i\gs 1}x^{\ast}(y_i)\va_i|>C^{-1}u).
 \ee
 This property is called weak tail domination. The weak tail domination can be understood 
 in terms of comparability of weak moments, i.e. for any integer $p\gs 1$ and $x^{\ast}\in B^{\ast}$
 \be\label{Ole2}
 \| \sum_{i\gs 1}x^{\ast}(x_i)\va_i\|_p \ls \bar{C}\| \sum_{i\gs 1}x^{\ast}(y_i)\va_i \|_p
 \ee  
 Oleszkiewicz asked whether or not it implies the comparability of strong moments.
 Namely whether (\ref{Ole1}) implies that
 \begin{align} 
 & \E \| \sum_{i\gs 1}x_i\va_i\|=\E \sup_{x^{\ast} \in B^{\ast}_1} \sum_{i\gs 1} x^{\ast}(x_i)\va_i \nonumber\\
 \label{Ole0}&\ls K\E \sup_{x^{\ast} \in B^{\ast}_1} \sum_{i\gs 1} x^{\ast}(y_i)\va_i= K \E \|\sum_{i\gs 1} y_i\va_i\|?
 \end{align}
 Note that in the Oleszkiewicz problem one may assume that $B$ is a separable space since we can easily restrict $B$ to the closure of $\mathbf{Lin}(y_1,x_1,y_2,x_2,\ldots)$. Therefore 
we have that
$$
\E \|\sum_{i\gs 1} y_i\va_i\|=\sup_{F\subset B^{\ast}_1}\E\sup_{x^{\ast}\in F}|\sum_{i\gs 1}x^{\ast}(y_i)\va_i|,
$$
where the supremum is taken over all finite sets $F$ contained in $B^{\ast}_1=\{x^{\ast}\in B^{\ast}:\;\|x^{\ast}\|\ls 1\}$.
We may assume that $\E \|\sum_{i\gs 1} y_i\va_i\|<\infty$ since otherwise there is nothing to prove.
Consequently for each $x^{\ast}\in B^{\ast}$ series $\sum_{i\gs 1}x^{\ast}(y_i)\va_i$ is convergent
which is equivalent to $\sum_{i\gs 1}(x^{\ast}(y_i))^2<\infty$.
Let $Q:B^{\ast}\ra \ell^2$ be defined by
$Q(x^{\ast})=(x^{\ast}(y_i))_{i\gs 1}$. It is clear that $Q:B^{\ast}/\ker Q\ra \ell^2$ is a linear
isomorphism on the closed subspace of $\ell^2$. First result which is an immediate consequence of
the Bernoulli theorem concerns the case when $Q$ maps $B^{\ast}$ onto $\ell^2$.
\bt\label{thm:2}
If $Q$ is onto $\ell^2$ then (\ref{Ole1}) implies (\ref{Ole0}).
\et
\begin{dwd}
Due to the Bernouli theorem i.e. (\ref{ole3}) there exists a decomposition of $Q(B^{\ast}_1)$,  
into $T_1,T_2$ such that $T_1+T_2\supset Q(B^{\ast}_1)$ and
\be\label{mela1}
b(Q(B^{\ast}_1))\gs L^{-1}(\sup_{t\in T_1}\|t\|_1+\gamma_G(T_2)).
\ee
Now for each $t\in T_1\cup T_2$ there exists a unique $[x^{\ast}_t]\in B^{\ast}/\ker Q$ such that
$Q([x^{\ast}_t])=t$. Let $R:B^{\ast}\ra \ell^2$ be defined $R(x^{\ast})=(x^{\ast}(x_i))_{i\gs 1}$.
Obviously (\ref{Ole1}) yields $\ker Q\subset \ker R$ and hence $R([x^{\ast}_t])$ is a unique point in   
$R(B^{\ast})$. Let $S_1=\{R([x^{\ast}_t])\in R(B^{\ast}):\;t\in T_1\}$ and
$S_2=\{R([x^{\ast}_t])\in R(B^{\ast}):\;t\in T_2\}$.
Since (\ref{Ole1}) or rather its consequence (\ref{Ole2}) implies
$$
\|R([x^{\ast}_t])-R([x^{\ast}_s])\|_2\ls \|Q([x^{\ast}_t])-Q([x^{\ast}_s])\|_2=\|t-s\|_2
$$
we get by Theorem \ref{thm:4} that
$$
\gamma_G(S_2)\ls K\gamma_G(T_2). 
$$
In the same way (\ref{Ole1}) gives that
$$
\sup_{s\in S_1}\|s\|_1=\sup_{t\in T_1}\|R([x^{\ast}_t])\|\ls K\sup_{t\in T_1}\|Q([x^{\ast}_t])\|_1=K\sup_{t\in T}\|t\|_1.
$$
Finally for each $x^{\ast}\in B^{\ast}_1$ there exits $t_1\in T_1$ and $t_2\in T_2$
such that $Q(x^{\ast})=t_1+t_2=Q([x^{\ast}_{t_1}])+Q([x^{\ast}_{t_2}])$.
Therefore $[x^{\ast}]=[x^{\ast}_{t_1}]+[x^{\ast}_{t_2}]$ and hence
$$
R([x^{\ast}])=R([x^{\ast}_{t_1}])+R([x^{\ast}_{t_2}]),
$$
which means that $S_1+S_2\supset R(B^{\ast})$. Since clearly by the easy upper bound part of the Bernoulli theorem \cite{Bed1}
$$
b(S_1)\ls L\sup_{s\in S_1}\|s\|_1,\;\;\mbox{and}\;\;b(S_2)\ls\gamma_G(S_2).
$$
we get by (\ref{mela1})
\begin{align*}
&b(R(B^{\ast}_1))\ls b(S_1)+b(S_2)\ls L(\sup_{s\in S_1}\|s\|_1+\gamma_G(S_2))\\
&\ls  KL(\sup_{t\in T_1}\|t\|_1+\gamma_G(T_2))\ls KL^2b(Q(B^{\ast}_1)).
\end{align*}
This ends the proof.
\end{dwd}
If $Q$ is not onto $\ell^2$ then the above argument fails but still it is believed that
the comparison holds. A partial result can be deduced from Theorem \ref{thm:1} namely
\begin{coro}
Suppose that for each $x^{\ast}\in B^{\ast}$ and $p\gs 1$
\be\label{ole6}
\|\sum_{i\gs 1} \min\{|x^{\ast}(x_i)|,1\}\va_i\|_{p}\ls C(p+\|\sum_{i\gs 1} \min\{|x^{\ast}(x_i)|,1\}\va_i\|_{p}).
\ee
Then (\ref{Ole0}) holds, i.e. 
$$
\E\|\sum_{i\gs 1}x_i\va_i\|\ls K\E\|\sum_{i\gs 1}y_i\va_i\|.
$$
\end{coro} 
\begin{dwd}
It suffices to notice that (\ref{ole6}) implies (\ref{ole2}) uand the apply Theorem \ref{thm:1}.
\end{dwd}

\end{document}